\documentclass[12pt]{artikel2}
\usepackage{a4wide,theorem,amstex,pstricks}

\let\phi\varphi
\newcommand{\ran}{\ensuremath{\mathop{\mathrm{Ran}}\nolimits}}
\newcommand{\unit}{\ensuremath{{\rm 1}\kern-.25em {\rm I}}}
\ifx\undefined\bysame
\newcommand{\bysame}{\leavevmode\hbox to3em{\hrulefill}\,}
\fi

\theoremstyle{plain}
\theoremheaderfont{\scshape}
\theorembodyfont{\slshape}
\newtheorem{theorem}{Theorem}
\newtheorem{lemma}[theorem]{Lemma}
\newtheorem{proposition}[theorem]{Proposition}

\newcommand{\QED}{\ifhmode\unskip\nobreak\fi\hfill
        \ifmmode\Box\else$\mathsurround=0pt\Box$\fi}
\newenvironment{proof}{\medbreak\parindent=0pt{\sc Proof:}\enspace}{%
        \nobreak\enspace\hbox{}\QED\par
        \ifdim\lastskip<\medskipamount\removelastskip\penalty55\medskip\fi
        \leftskip=0pt}

\title{An obstruction for {\boldmath $q$-deformation} of the
       convolution product}
\author{Hans van Leeuwen and Hans Maassen}
\renewcommand{\today}{{\footnotesize 
        Department of mathematics KUN\\
        Toernooiveld 1\\
        6525~ED~~Nijmegen\\
        the Netherlands\\
        E-mail: leeuwen@@sci.kun.nl, maassen@@sci.kun.nl}}
\begin{document}
\maketitle
\begin{abstract}
We consider two independent $q$-Gaussian random variables $X_0$ and $X_1$
and a function $\gamma$ chosen in such a way that $\gamma(X_0)$ and
$X_0$ have the same distribution. For $q \in (0, 1)$ we find that at 
least the fourth moments of $X_0 + X_1$ and $\gamma(X_0) + X_1$
are different. We conclude that no $q$-deformed convolution
product can exist for functions of independent $q$-Gaussian random
variables.\par
\smallskip
{\footnotesize\noindent
1995 PACS numbers: 02.50.Cw, 05.40.+j, 03.65.Db, 42.50.Lc\par\noindent
}
\end{abstract}

\section{Introduction and Notation}

In 1991 Bo\.{z}ejko and Speicher introduced a deformation of Brownian
motion by a parameter $q \in [-1, 1]$ (cf.~\cite{bs91,bs92}). Their
construction is based on a $q$-deformation, $\mathcal{F}_q(\mathcal{H})$, of
the full Fock space over a separable Hilbert space $\mathcal{H}$.  Their
random variables are given by self-adjoint operators of the form
\[
X(f):= a(f) + a(f)^*,\qquad f\in\mathcal{H},
\]
where $a(f)$ and $a(f)^*$ are the annihilation and creation operators
associated to $f$ satisfying the $q$-deformed commutation relation,
\begin{equation}
\label{comm1}
a(f) a(g)^* - q a(g)^* a(f) = \langle f, g \rangle \unit.
\end{equation}
This commutation relation was first introduced by Frisch and Bourret
in~\cite{fb70} and various aspects of it are studied in
a.o.~\cite{dn93,gre91,jsw94,lm95}.

The above construction was for some time considered a good candidate
for a $q$-deformed notion of the concept of independence
itself. Indeed for $q=1$ the random variables $X(f)$ and $X(g)$ with
$f \perp g$ are independent Gaussian random variables in the classical
sense, for $q=0$ they are freely independent in the sense of
Voiculescu~\cite{vdn92}.  In both cases a convolution law holds for
sums of functions of $X(f)$ and $X(g)$. For $q=1$ the convolution is
ordinary convolution whereas for $q=0$ the convolution is found to be
an interesting operation involving Gauchy transforms and inverted
functions~\cite{maa92,vdn92}.  In this paper we show, by a simple
example, that for $q \in (0, 1)$ no such convolution law holds since
the distributions of functions of $X(f)$ and $X(g)$ do not determine
the distribution of their sum.

The construction of the Fock representation for (\ref{comm1}) is
described in~\cite{bs91,dn93} but for completeness we give the
necessary details here. Operators $a(f)$ and $a(f)^*$ are, for all $f
\in \mathcal{H}$, defined on the full Fock space $\mathcal{F} := \mathbb{C}
\oplus \bigoplus_{n=1}^\infty \mathcal{H}^{\otimes n}$ by:
\[
a(f)^* h_1 \otimes \cdots \otimes h_n
:= 
        f \otimes h_1 \otimes \cdots \otimes h_n,\qquad n\in \mathbb{N},
        h_1, \ldots, h_n \in \mathcal{H}
\]
and
\begin{align}
\label{vacuum_eis}
a(f) \Omega
& :=
        0,\\
a(f) h_1 \otimes \cdots \otimes h_n
& :=
        \sum_{k=1}^n q^{k-1} \langle f, h_k \rangle
        h_1 \otimes \cdots \check{h}_k \cdots \otimes h_n,
        \qquad n \geq 1\nonumber
\end{align}
where the notation $h_1 \otimes \cdots \check{h}_k \cdots \otimes h_n$
stands for the tensor product $h_1 \otimes \cdots \otimes h_{k-1}
\otimes h_{k+1} \otimes \cdots \otimes h_n$ and $\Omega = 1 \oplus 0
\oplus 0 \oplus \cdots$.  In order to ensure that $a(f)^*$ is the
adjoint of $a(f)$ for all $f \in \mathcal{H}$, Bo\.{z}ejko and Speicher
recursively define an inner product $\langle \cdot, \cdot\rangle_q$ on
$\mathcal{F}$ as:
\begin{multline*}
\langle g_1 \otimes \cdots \otimes g_m, h_1 \otimes \cdots \otimes h_n\rangle_q
:=
        \delta_{n, m}\langle g_2 \otimes \cdots \otimes g_m,
        a(g_1) h_1 \otimes \cdots \otimes h_n\rangle_q\\
=
        \delta_{n, m}\sum_{k=1}^n q^{k-1} \langle g_1, h_k \rangle
        \langle g_2 \otimes \cdots \otimes g_m,
        h_1 \otimes \cdots \check{h}_k \cdots \otimes h_n \rangle_q.
\end{multline*}
We denote the full Fock space $\mathcal{F}$ equipped with this inner
product by $\mathcal{F}_q(\mathcal{H})$.  By the GNS construction
there exists, up to unitary equivalence, only one cyclic
representation of the relations (\ref{comm1}) and
(\ref{vacuum_eis}). For $\mathcal{H} = \mathbb{C}$ the above
construction reduces to $\mathcal{F}_q(\mathbb{C}) \cong l^2(\mathbb{N},
[n]_q!)$, where $[n]_q = (1-q^n)/(1-q)$ and $[n]_q! = \prod_{j=1}^n
[j]_q$ with $[0]_q! = 1$.

In~\cite{bs92,lm95} the density of the $q$-Gaussian distribution,
$\nu_q(dx)$, of the random variable $X_0 = a(f_0) + a(f_0)^*$ with $f_0 \in
\mathcal{H}$ and $\| f_0 \| = 1$ is calculated. This density is a measure
on $\mathbb{R}$, where it is supported on the interval $[-2/\sqrt{1-q},
2/\sqrt{1-q}]$. If we denote the $n$-fold product $\prod_{k=0}^{n-1}
(1 - a q^k)$ by $(a; q)_n$ and agree on $(a_1, \ldots, a_m; q)_n =
(a_1; q)_n \cdots (a_m; q)_n$, then $\nu_q(dx)$ can be written as:
\[
\nu_q(dx)
=
	\nu_q^\prime(x) dx
=
        \frac{1}{\pi} \sqrt{1-q} \sin{\theta} 
        (q, q v^2, q v^{-2}; q)_\infty \, dx,
\]
where $2\cos{\theta} = x\sqrt{1-q}$ and $v = \exp{(i\theta)}$.

To state the main theorem of this paper we define $X_1$ to be the
random variable $a(f_1) + a(f_1)^*$ for some $f_1 \in \mathcal{H}$ with
$\| f_1 \| = 1$ and $\langle f_0, f_1\rangle = 0$. Then $X_0$ and $X_1$ are
$q$-Gaussian random variables, independent in the sense of quantum
probability (\cite{bs91,kum}).
\begin{theorem}
\label{hoofdstelling}
There exists a function $\gamma\colon \mathbb{R} \rightarrow \mathbb{R}$
such that $X_0$ and $\gamma(X_0)$ are identically distributed but $X_0+X_1$
and $\gamma(X_0)+X_1$ are not.
\end{theorem}
The consequence of this theorem is that the distribution of the sum of
two or more random variables depends on the choice of random variables
and not solely on the respective distributions of these random
variables. This means that a $q$-convolution paralleling the known
convolution for probability measures for the cases $q=0$
(cf.~\cite{maa92,vdn92}) and $q = 1$ cannot exist.

In contrast to the above, Nica, in~\cite{ni95}, constructs a convolution law
for probability distributions that interpolates between the known
cases $q=0$ and $q=1$. Theorem~\ref{hoofdstelling} implies that
this interpolation is incompatible with relation
(\ref{comm1}). In fact this can also be seen by explicit calculation of the
moments of the distribution of $X_0^n + X_1^m$, $n, m \in \mathbb{N}$ using
the convolution law Nica suggests and using the structure inherently
present in $\mathcal{F}_q(\mathcal{H})$. From the fourth moment onwards the
moments differ for $n, m \geq 1$, although they are the same for the
cases $q = 0$ and $q = \pm 1$, as they should be.

In the next section we will prove theorem~\ref{hoofdstelling} by
constructing the function $\gamma$ and showing that the fourth moment
of the distribution of $\gamma(X_0)+X_1$ is strictly smaller then the
fourth moment of the distribution of $X_0 + X_1$ for $q \in (0, 1)$.

\section{Construction of {\boldmath $\gamma$} and proof of theorem}

In \cite{lm95} we construct the unitary operator
$U\colon \mathcal{F}_q(\mathbb{C}) \rightarrow L^2(\mathbb{R}, \nu_q)$
that diagonalizes the operator $X = a + a^*$ with $a = a(1)$,
such that $U X = T U$ with $T$ the operator of pointwise
multiplication on $L^2(\mathbb{R}, \nu_q)$ given by $(T f)(x) = x f(x)$
for $f \in L^2(\mathbb{R}, \nu_q)$.

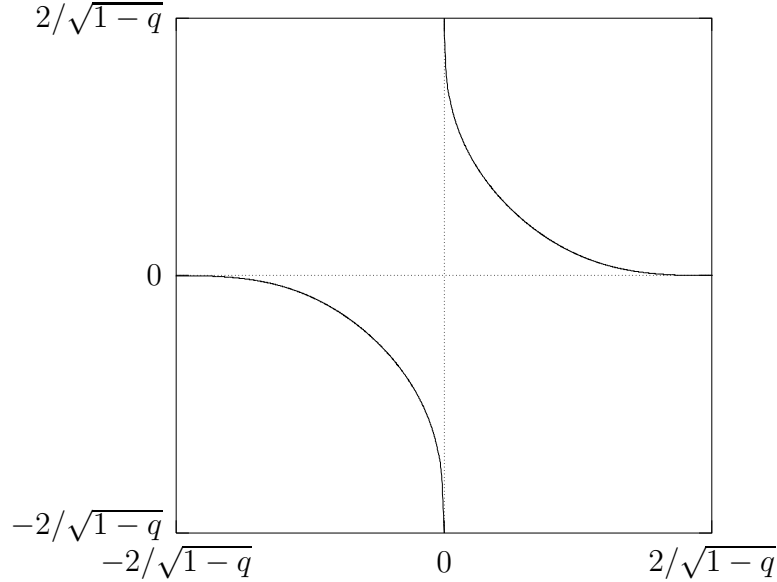
\begin{figure}
\centerline{
\psset{unit=0.240900pt}
\psset{arrowsize=3pt 3.2 1.4 .3}
\begin{pspicture}(0,0)(1081,900)
\psset{linewidth=0.3pt,linestyle=dotted,dotsep=1pt}\psline(176,473)(1017,473)
\psset{linewidth=0.3pt,linestyle=dotted,dotsep=1pt}\psline(597,68)(597,877)
\psset{linewidth=0.35pt,linestyle=solid}\psline(176,68)(196,68)
\psset{linewidth=0.35pt,linestyle=solid}\psline(1017,68)(997,68)
\rput[r](154,80){$-2/\sqrt{1-q}$}
\psset{linewidth=0.35pt,linestyle=solid}\psline(176,473)(196,473)
\psset{linewidth=0.35pt,linestyle=solid}\psline(1017,473)(997,473)
\rput[r](154,473){$0$}
\psset{linewidth=0.35pt,linestyle=solid}\psline(176,877)(196,877)
\psset{linewidth=0.35pt,linestyle=solid}\psline(1017,877)(997,877)
\rput[r](154,877){$2/\sqrt{1-q}$}
\psset{linewidth=0.35pt,linestyle=solid}\psline(176,68)(176,88)
\psset{linewidth=0.35pt,linestyle=solid}\psline(176,877)(176,857)
\rput(176,23){$-2/\sqrt{1-q}$}
\psset{linewidth=0.35pt,linestyle=solid}\psline(597,68)(597,88)
\psset{linewidth=0.35pt,linestyle=solid}\psline(597,877)(597,857)
\rput(597,23){$0$}
\psset{linewidth=0.35pt,linestyle=solid}\psline(1017,68)(1017,88)
\psset{linewidth=0.35pt,linestyle=solid}\psline(1017,877)(1017,857)
\rput(1017,23){$2/\sqrt{1-q}$}
\psset{linewidth=0.35pt,linestyle=solid}\psline(176,68)(1017,68)(1017,877)(176,877)(176,68)
\psset{linewidth=0.35pt,linestyle=solid}\psline(176,473)(176,473)(176,472)(177,472)(177,472)(178,472)(178,472)(179,472)(179,472)(179,472)(180,472)(180,472)(181,472)(181,472)(181,472)(182,472)(182,472)(183,472)(183,472)(184,472)(184,472)(184,472)(185,472)(185,472)(186,472)(186,472)(187,472)(187,472)(187,472)(188,472)(188,472)(189,472)(189,472)(189,472)(190,472)(190,472)(191,472)(191,472)(192,472)(192,472)(192,472)(193,472)(193,472)(194,472)(194,472)(195,472)(195,472)(195,472)(196,472)(196,472)(197,472)
\psset{linewidth=0.35pt,linestyle=solid}\psline(197,472)(197,472)(197,472)(198,472)(198,472)(199,472)(199,472)(200,472)(200,472)(200,472)(201,472)(201,472)(202,472)(202,472)(202,472)(203,472)(203,472)(204,472)(204,472)(205,472)(205,472)(205,472)(206,472)(206,472)(207,472)(207,472)(208,472)(208,472)(208,472)(209,472)(209,472)(210,472)(210,472)(210,472)(211,472)(211,472)(212,472)(212,472)(213,472)(213,472)(213,472)(214,472)(214,472)(215,472)(215,472)(216,472)(216,472)(216,472)(217,472)(217,472)(218,472)
\psset{linewidth=0.35pt,linestyle=solid}\psline(218,472)(218,472)(218,472)(219,472)(219,472)(220,472)(220,472)(221,472)(221,472)(221,472)(222,472)(222,472)(223,472)(223,472)(224,472)(224,472)(224,472)(225,472)(225,472)(226,472)(226,472)(226,472)(227,472)(227,472)(228,472)(228,472)(229,472)(229,472)(229,472)(230,472)(230,472)(231,472)(231,472)(232,472)(232,472)(232,472)(233,472)(233,472)(234,472)(234,471)(234,471)(235,471)(235,471)(236,471)(236,471)(237,471)(237,471)(237,471)(238,471)(238,471)(239,471)
\psset{linewidth=0.35pt,linestyle=solid}\psline(239,471)(239,471)(239,471)(240,471)(240,471)(241,471)(241,471)(242,471)(242,471)(242,471)(243,471)(243,471)(244,471)(244,471)(245,471)(245,471)(245,471)(246,471)(246,471)(247,471)(247,471)(247,471)(248,471)(248,471)(249,471)(249,471)(250,471)(250,471)(250,471)(251,471)(251,471)(252,471)(252,470)(253,470)(253,470)(253,470)(254,470)(254,470)(255,470)(255,470)(255,470)(256,470)(256,470)(257,470)(257,470)(258,470)(258,470)(258,470)(259,470)(259,470)(260,470)
\psset{linewidth=0.35pt,linestyle=solid}\psline(260,470)(260,470)(261,470)(261,470)(261,470)(262,470)(262,470)(263,470)(263,470)(263,470)(264,470)(264,470)(265,470)(265,469)(266,469)(266,469)(266,469)(267,469)(267,469)(268,469)(268,469)(269,469)(269,469)(269,469)(270,469)(270,469)(271,469)(271,469)(271,469)(272,469)(272,469)(273,469)(273,469)(274,469)(274,469)(274,469)(275,469)(275,468)(276,468)(276,468)(276,468)(277,468)(277,468)(278,468)(278,468)(279,468)(279,468)(279,468)(280,468)(280,468)(281,468)
\psset{linewidth=0.35pt,linestyle=solid}\psline(281,468)(281,468)(282,468)(282,468)(282,468)(283,468)(283,468)(284,467)(284,467)(284,467)(285,467)(285,467)(286,467)(286,467)(287,467)(287,467)(287,467)(288,467)(288,467)(289,467)(289,467)(290,467)(290,467)(290,467)(291,467)(291,466)(292,466)(292,466)(292,466)(293,466)(293,466)(294,466)(294,466)(295,466)(295,466)(295,466)(296,466)(296,466)(297,466)(297,466)(298,466)(298,465)(298,465)(299,465)(299,465)(300,465)(300,465)(300,465)(301,465)(301,465)(302,465)
\psset{linewidth=0.35pt,linestyle=solid}\psline(302,465)(302,465)(303,465)(303,465)(303,465)(304,464)(304,464)(305,464)(305,464)(306,464)(306,464)(306,464)(307,464)(307,464)(308,464)(308,464)(308,464)(309,464)(309,463)(310,463)(310,463)(311,463)(311,463)(311,463)(312,463)(312,463)(313,463)(313,463)(314,463)(314,463)(314,462)(315,462)(315,462)(316,462)(316,462)(316,462)(317,462)(317,462)(318,462)(318,462)(319,462)(319,462)(319,461)(320,461)(320,461)(321,461)(321,461)(321,461)(322,461)(322,461)(323,461)
\psset{linewidth=0.35pt,linestyle=solid}\psline(323,461)(323,461)(324,461)(324,460)(324,460)(325,460)(325,460)(326,460)(326,460)(327,460)(327,460)(327,460)(328,460)(328,459)(329,459)(329,459)(329,459)(330,459)(330,459)(331,459)(331,459)(332,459)(332,458)(332,458)(333,458)(333,458)(334,458)(334,458)(335,458)(335,458)(335,458)(336,457)(336,457)(337,457)(337,457)(337,457)(338,457)(338,457)(339,457)(339,457)(340,456)(340,456)(340,456)(341,456)(341,456)(342,456)(342,456)(343,456)(343,456)(343,455)(344,455)
\psset{linewidth=0.35pt,linestyle=solid}\psline(344,455)(344,455)(345,455)(345,455)(345,455)(346,455)(346,455)(347,454)(347,454)(348,454)(348,454)(348,454)(349,454)(349,454)(350,453)(350,453)(351,453)(351,453)(351,453)(352,453)(352,453)(353,453)(353,452)(353,452)(354,452)(354,452)(355,452)(355,452)(356,452)(356,451)(356,451)(357,451)(357,451)(358,451)(358,451)(358,451)(359,450)(359,450)(360,450)(360,450)(361,450)(361,450)(361,450)(362,449)(362,449)(363,449)(363,449)(364,449)(364,449)(364,449)(365,448)
\psset{linewidth=0.35pt,linestyle=solid}\psline(365,448)(365,448)(366,448)(366,448)(366,448)(367,448)(367,447)(368,447)(368,447)(369,447)(369,447)(369,447)(370,447)(370,446)(371,446)(371,446)(372,446)(372,446)(372,446)(373,445)(373,445)(374,445)(374,445)(374,445)(375,445)(375,444)(376,444)(376,444)(377,444)(377,444)(377,444)(378,443)(378,443)(379,443)(379,443)(380,443)(380,442)(380,442)(381,442)(381,442)(382,442)(382,442)(382,441)(383,441)(383,441)(384,441)(384,441)(385,440)(385,440)(385,440)(386,440)
\psset{linewidth=0.35pt,linestyle=solid}\psline(386,440)(386,440)(387,440)(387,439)(388,439)(388,439)(388,439)(389,439)(389,438)(390,438)(390,438)(390,438)(391,438)(391,437)(392,437)(392,437)(393,437)(393,437)(393,436)(394,436)(394,436)(395,436)(395,436)(396,435)(396,435)(396,435)(397,435)(397,435)(398,434)(398,434)(398,434)(399,434)(399,434)(400,433)(400,433)(401,433)(401,433)(401,433)(402,432)(402,432)(403,432)(403,432)(403,431)(404,431)(404,431)(405,431)(405,431)(406,430)(406,430)(406,430)(407,430)
\psset{linewidth=0.35pt,linestyle=solid}\psline(407,430)(407,430)(408,429)(408,429)(409,429)(409,429)(409,428)(410,428)(410,428)(411,428)(411,427)(411,427)(412,427)(412,427)(413,427)(413,426)(414,426)(414,426)(414,426)(415,425)(415,425)(416,425)(416,425)(417,424)(417,424)(417,424)(418,424)(418,423)(419,423)(419,423)(419,423)(420,422)(420,422)(421,422)(421,422)(422,421)(422,421)(422,421)(423,421)(423,420)(424,420)(424,420)(425,420)(425,419)(425,419)(426,419)(426,419)(427,418)(427,418)(427,418)(428,418)
\psset{linewidth=0.35pt,linestyle=solid}\psline(428,418)(428,417)(429,417)(429,417)(430,417)(430,416)(430,416)(431,416)(431,416)(432,415)(432,415)(433,415)(433,414)(433,414)(434,414)(434,414)(435,413)(435,413)(435,413)(436,413)(436,412)(437,412)(437,412)(438,411)(438,411)(438,411)(439,411)(439,410)(440,410)(440,410)(440,409)(441,409)(441,409)(442,409)(442,408)(443,408)(443,408)(443,407)(444,407)(444,407)(445,407)(445,406)(446,406)(446,406)(446,405)(447,405)(447,405)(448,404)(448,404)(448,404)(449,404)
\psset{linewidth=0.35pt,linestyle=solid}\psline(449,404)(449,403)(450,403)(450,403)(451,402)(451,402)(451,402)(452,401)(452,401)(453,401)(453,400)(454,400)(454,400)(454,400)(455,399)(455,399)(456,399)(456,398)(456,398)(457,398)(457,397)(458,397)(458,397)(459,396)(459,396)(459,396)(460,395)(460,395)(461,395)(461,394)(462,394)(462,394)(462,393)(463,393)(463,393)(464,392)(464,392)(464,392)(465,391)(465,391)(466,391)(466,390)(467,390)(467,390)(467,389)(468,389)(468,389)(469,388)(469,388)(470,388)(470,387)
\psset{linewidth=0.35pt,linestyle=solid}\psline(470,387)(470,387)(471,386)(471,386)(472,386)(472,385)(472,385)(473,385)(473,384)(474,384)(474,384)(475,383)(475,383)(475,383)(476,382)(476,382)(477,381)(477,381)(477,381)(478,380)(478,380)(479,380)(479,379)(480,379)(480,378)(480,378)(481,378)(481,377)(482,377)(482,377)(483,376)(483,376)(483,375)(484,375)(484,375)(485,374)(485,374)(485,373)(486,373)(486,373)(487,372)(487,372)(488,371)(488,371)(488,371)(489,370)(489,370)(490,369)(490,369)(491,369)(491,368)
\psset{linewidth=0.35pt,linestyle=solid}\psline(491,368)(491,368)(492,367)(492,367)(493,367)(493,366)(493,366)(494,365)(494,365)(495,364)(495,364)(496,364)(496,363)(496,363)(497,362)(497,362)(498,362)(498,361)(499,361)(499,360)(499,360)(500,359)(500,359)(501,358)(501,358)(501,358)(502,357)(502,357)(503,356)(503,356)(504,355)(504,355)(504,354)(505,354)(505,354)(506,353)(506,353)(507,352)(507,352)(507,351)(508,351)(508,350)(509,350)(509,349)(509,349)(510,348)(510,348)(511,348)(511,347)(512,347)(512,346)
\psset{linewidth=0.35pt,linestyle=solid}\psline(512,346)(512,346)(513,345)(513,345)(514,344)(514,344)(515,343)(515,343)(515,342)(516,342)(516,341)(517,341)(517,340)(517,340)(518,339)(518,339)(519,338)(519,338)(520,337)(520,337)(520,336)(521,336)(521,335)(522,335)(522,334)(522,333)(523,333)(523,332)(524,332)(524,331)(525,331)(525,330)(525,330)(526,329)(526,329)(527,328)(527,328)(528,327)(528,326)(528,326)(529,325)(529,325)(530,324)(530,324)(530,323)(531,323)(531,322)(532,321)(532,321)(533,320)(533,320)
\psset{linewidth=0.35pt,linestyle=solid}\psline(533,320)(533,319)(534,318)(534,318)(535,317)(535,317)(536,316)(536,315)(536,315)(537,314)(537,314)(538,313)(538,312)(538,312)(539,311)(539,311)(540,310)(540,309)(541,309)(541,308)(541,307)(542,307)(542,306)(543,305)(543,305)(544,304)(544,303)(544,303)(545,302)(545,301)(546,301)(546,300)(546,299)(547,299)(547,298)(548,297)(548,297)(549,296)(549,295)(549,295)(550,294)(550,293)(551,292)(551,292)(552,291)(552,290)(552,290)(553,289)(553,288)(554,287)(554,287)
\psset{linewidth=0.35pt,linestyle=solid}\psline(554,287)(554,286)(555,285)(555,284)(556,283)(556,283)(557,282)(557,281)(557,280)(558,279)(558,279)(559,278)(559,277)(559,276)(560,275)(560,275)(561,274)(561,273)(562,272)(562,271)(562,270)(563,269)(563,268)(564,268)(564,267)(565,266)(565,265)(565,264)(566,263)(566,262)(567,261)(567,260)(567,259)(568,258)(568,257)(569,256)(569,255)(570,254)(570,253)(570,252)(571,251)(571,250)(572,249)(572,248)(573,247)(573,246)(573,245)(574,243)(574,242)(575,241)(575,240)
\psset{linewidth=0.35pt,linestyle=solid}\psline(575,240)(575,239)(576,238)(576,236)(577,235)(577,234)(578,233)(578,231)(578,230)(579,229)(579,227)(580,226)(580,224)(581,223)(581,222)(581,220)(582,219)(582,217)(583,216)(583,214)(583,212)(584,211)(584,209)(585,207)(585,205)(586,203)(586,201)(586,199)(587,197)(587,195)(588,193)(588,191)(589,190)(589,188)(589,186)(590,184)(590,181)(591,179)(591,176)(591,172)(592,168)(592,164)(593,158)(593,152)(594,145)(594,138)(594,129)(595,119)(595,108)(596,96)(596,83)
\psset{linewidth=0.35pt,linestyle=solid}\psline(596,83)(597,68)
\psset{linewidth=0.35pt,linestyle=solid}\psline(597,877)(597,877)(597,862)(597,849)(598,837)(598,826)(599,816)(599,807)(599,800)(600,793)(600,787)(601,781)(601,777)(602,773)(602,769)(602,766)(603,764)(603,761)(604,759)(604,757)(604,755)(605,754)(605,752)(606,750)(606,748)(607,746)(607,744)(607,742)(608,740)(608,738)(609,736)(609,734)(610,733)(610,731)(610,729)(611,728)(611,726)(612,725)(612,723)(612,722)(613,721)(613,719)(614,718)(614,716)(615,715)(615,714)(615,712)(616,711)(616,710)(617,709)(617,707)
\psset{linewidth=0.35pt,linestyle=solid}\psline(617,707)(618,706)(618,705)(618,704)(619,703)(619,702)(620,700)(620,699)(620,698)(621,697)(621,696)(622,695)(622,694)(623,693)(623,692)(623,691)(624,690)(624,689)(625,688)(625,687)(626,686)(626,685)(626,684)(627,683)(627,682)(628,681)(628,680)(628,679)(629,678)(629,677)(630,677)(630,676)(631,675)(631,674)(631,673)(632,672)(632,671)(633,670)(633,670)(634,669)(634,668)(634,667)(635,666)(635,666)(636,665)(636,664)(636,663)(637,662)(637,662)(638,661)(638,660)
\psset{linewidth=0.35pt,linestyle=solid}\psline(638,660)(639,659)(639,658)(639,658)(640,657)(640,656)(641,655)(641,655)(641,654)(642,653)(642,653)(643,652)(643,651)(644,650)(644,650)(644,649)(645,648)(645,648)(646,647)(646,646)(647,646)(647,645)(647,644)(648,644)(648,643)(649,642)(649,642)(649,641)(650,640)(650,640)(651,639)(651,638)(652,638)(652,637)(652,636)(653,636)(653,635)(654,634)(654,634)(655,633)(655,633)(655,632)(656,631)(656,631)(657,630)(657,630)(657,629)(658,628)(658,628)(659,627)(659,627)
\psset{linewidth=0.35pt,linestyle=solid}\psline(659,627)(660,626)(660,625)(660,625)(661,624)(661,624)(662,623)(662,622)(663,622)(663,621)(663,621)(664,620)(664,620)(665,619)(665,619)(665,618)(666,617)(666,617)(667,616)(667,616)(668,615)(668,615)(668,614)(669,614)(669,613)(670,613)(670,612)(671,612)(671,611)(671,610)(672,610)(672,609)(673,609)(673,608)(673,608)(674,607)(674,607)(675,606)(675,606)(676,605)(676,605)(676,604)(677,604)(677,603)(678,603)(678,602)(678,602)(679,601)(679,601)(680,600)(680,600)
\psset{linewidth=0.35pt,linestyle=solid}\psline(680,600)(681,599)(681,599)(681,598)(682,598)(682,597)(683,597)(683,597)(684,596)(684,596)(684,595)(685,595)(685,594)(686,594)(686,593)(686,593)(687,592)(687,592)(688,591)(688,591)(689,591)(689,590)(689,590)(690,589)(690,589)(691,588)(691,588)(692,587)(692,587)(692,587)(693,586)(693,586)(694,585)(694,585)(694,584)(695,584)(695,583)(696,583)(696,583)(697,582)(697,582)(697,581)(698,581)(698,581)(699,580)(699,580)(700,579)(700,579)(700,578)(701,578)(701,578)
\psset{linewidth=0.35pt,linestyle=solid}\psline(701,578)(702,577)(702,577)(702,576)(703,576)(703,576)(704,575)(704,575)(705,574)(705,574)(705,574)(706,573)(706,573)(707,572)(707,572)(708,572)(708,571)(708,571)(709,570)(709,570)(710,570)(710,569)(710,569)(711,568)(711,568)(712,568)(712,567)(713,567)(713,567)(713,566)(714,566)(714,565)(715,565)(715,565)(716,564)(716,564)(716,564)(717,563)(717,563)(718,562)(718,562)(718,562)(719,561)(719,561)(720,561)(720,560)(721,560)(721,560)(721,559)(722,559)(722,559)
\psset{linewidth=0.35pt,linestyle=solid}\psline(722,559)(723,558)(723,558)(723,557)(724,557)(724,557)(725,556)(725,556)(726,556)(726,555)(726,555)(727,555)(727,554)(728,554)(728,554)(729,553)(729,553)(729,553)(730,552)(730,552)(731,552)(731,551)(731,551)(732,551)(732,550)(733,550)(733,550)(734,549)(734,549)(734,549)(735,548)(735,548)(736,548)(736,547)(737,547)(737,547)(737,546)(738,546)(738,546)(739,545)(739,545)(739,545)(740,545)(740,544)(741,544)(741,544)(742,543)(742,543)(742,543)(743,542)(743,542)
\psset{linewidth=0.35pt,linestyle=solid}\psline(743,542)(744,542)(744,541)(745,541)(745,541)(745,541)(746,540)(746,540)(747,540)(747,539)(747,539)(748,539)(748,538)(749,538)(749,538)(750,538)(750,537)(750,537)(751,537)(751,536)(752,536)(752,536)(753,536)(753,535)(753,535)(754,535)(754,534)(755,534)(755,534)(755,534)(756,533)(756,533)(757,533)(757,532)(758,532)(758,532)(758,532)(759,531)(759,531)(760,531)(760,531)(760,530)(761,530)(761,530)(762,529)(762,529)(763,529)(763,529)(763,528)(764,528)(764,528)
\psset{linewidth=0.35pt,linestyle=solid}\psline(764,528)(765,528)(765,527)(766,527)(766,527)(766,527)(767,526)(767,526)(768,526)(768,526)(768,525)(769,525)(769,525)(770,525)(770,524)(771,524)(771,524)(771,524)(772,523)(772,523)(773,523)(773,523)(774,522)(774,522)(774,522)(775,522)(775,521)(776,521)(776,521)(776,521)(777,520)(777,520)(778,520)(778,520)(779,519)(779,519)(779,519)(780,519)(780,518)(781,518)(781,518)(782,518)(782,518)(782,517)(783,517)(783,517)(784,517)(784,516)(784,516)(785,516)(785,516)
\psset{linewidth=0.35pt,linestyle=solid}\psline(785,516)(786,515)(786,515)(787,515)(787,515)(787,515)(788,514)(788,514)(789,514)(789,514)(790,514)(790,513)(790,513)(791,513)(791,513)(792,512)(792,512)(792,512)(793,512)(793,512)(794,511)(794,511)(795,511)(795,511)(795,511)(796,510)(796,510)(797,510)(797,510)(797,510)(798,509)(798,509)(799,509)(799,509)(800,509)(800,508)(800,508)(801,508)(801,508)(802,508)(802,507)(803,507)(803,507)(803,507)(804,507)(804,506)(805,506)(805,506)(805,506)(806,506)(806,505)
\psset{linewidth=0.35pt,linestyle=solid}\psline(806,505)(807,505)(807,505)(808,505)(808,505)(808,505)(809,504)(809,504)(810,504)(810,504)(811,504)(811,503)(811,503)(812,503)(812,503)(813,503)(813,503)(813,502)(814,502)(814,502)(815,502)(815,502)(816,501)(816,501)(816,501)(817,501)(817,501)(818,501)(818,500)(819,500)(819,500)(819,500)(820,500)(820,500)(821,499)(821,499)(821,499)(822,499)(822,499)(823,499)(823,498)(824,498)(824,498)(824,498)(825,498)(825,498)(826,498)(826,497)(827,497)(827,497)(827,497)
\psset{linewidth=0.35pt,linestyle=solid}\psline(827,497)(828,497)(828,497)(829,496)(829,496)(829,496)(830,496)(830,496)(831,496)(831,496)(832,495)(832,495)(832,495)(833,495)(833,495)(834,495)(834,495)(835,494)(835,494)(835,494)(836,494)(836,494)(837,494)(837,494)(837,493)(838,493)(838,493)(839,493)(839,493)(840,493)(840,493)(840,492)(841,492)(841,492)(842,492)(842,492)(842,492)(843,492)(843,492)(844,491)(844,491)(845,491)(845,491)(845,491)(846,491)(846,491)(847,490)(847,490)(848,490)(848,490)(848,490)
\psset{linewidth=0.35pt,linestyle=solid}\psline(848,490)(849,490)(849,490)(850,490)(850,489)(850,489)(851,489)(851,489)(852,489)(852,489)(853,489)(853,489)(853,489)(854,488)(854,488)(855,488)(855,488)(856,488)(856,488)(856,488)(857,488)(857,488)(858,487)(858,487)(858,487)(859,487)(859,487)(860,487)(860,487)(861,487)(861,487)(861,486)(862,486)(862,486)(863,486)(863,486)(864,486)(864,486)(864,486)(865,486)(865,485)(866,485)(866,485)(866,485)(867,485)(867,485)(868,485)(868,485)(869,485)(869,485)(869,484)
\psset{linewidth=0.35pt,linestyle=solid}\psline(869,484)(870,484)(870,484)(871,484)(871,484)(872,484)(872,484)(872,484)(873,484)(873,484)(874,484)(874,483)(874,483)(875,483)(875,483)(876,483)(876,483)(877,483)(877,483)(877,483)(878,483)(878,483)(879,483)(879,482)(879,482)(880,482)(880,482)(881,482)(881,482)(882,482)(882,482)(882,482)(883,482)(883,482)(884,482)(884,481)(885,481)(885,481)(885,481)(886,481)(886,481)(887,481)(887,481)(887,481)(888,481)(888,481)(889,481)(889,481)(890,480)(890,480)(890,480)
\psset{linewidth=0.35pt,linestyle=solid}\psline(890,480)(891,480)(891,480)(892,480)(892,480)(893,480)(893,480)(893,480)(894,480)(894,480)(895,480)(895,480)(895,479)(896,479)(896,479)(897,479)(897,479)(898,479)(898,479)(898,479)(899,479)(899,479)(900,479)(900,479)(901,479)(901,479)(901,479)(902,479)(902,478)(903,478)(903,478)(903,478)(904,478)(904,478)(905,478)(905,478)(906,478)(906,478)(906,478)(907,478)(907,478)(908,478)(908,478)(909,478)(909,478)(909,478)(910,477)(910,477)(911,477)(911,477)(911,477)
\psset{linewidth=0.35pt,linestyle=solid}\psline(911,477)(912,477)(912,477)(913,477)(913,477)(914,477)(914,477)(914,477)(915,477)(915,477)(916,477)(916,477)(917,477)(917,477)(917,477)(918,477)(918,476)(919,476)(919,476)(919,476)(920,476)(920,476)(921,476)(921,476)(922,476)(922,476)(922,476)(923,476)(923,476)(924,476)(924,476)(924,476)(925,476)(925,476)(926,476)(926,476)(927,476)(927,476)(927,476)(928,476)(928,475)(929,475)(929,475)(930,475)(930,475)(930,475)(931,475)(931,475)(932,475)(932,475)(932,475)
\psset{linewidth=0.35pt,linestyle=solid}\psline(932,475)(933,475)(933,475)(934,475)(934,475)(935,475)(935,475)(935,475)(936,475)(936,475)(937,475)(937,475)(938,475)(938,475)(938,475)(939,475)(939,475)(940,475)(940,475)(940,475)(941,475)(941,474)(942,474)(942,474)(943,474)(943,474)(943,474)(944,474)(944,474)(945,474)(945,474)(946,474)(946,474)(946,474)(947,474)(947,474)(948,474)(948,474)(948,474)(949,474)(949,474)(950,474)(950,474)(951,474)(951,474)(951,474)(952,474)(952,474)(953,474)(953,474)(954,474)
\psset{linewidth=0.35pt,linestyle=solid}\psline(954,474)(954,474)(954,474)(955,474)(955,474)(956,474)(956,474)(956,474)(957,474)(957,474)(958,474)(958,474)(959,474)(959,474)(959,473)(960,473)(960,473)(961,473)(961,473)(961,473)(962,473)(962,473)(963,473)(963,473)(964,473)(964,473)(964,473)(965,473)(965,473)(966,473)(966,473)(967,473)(967,473)(967,473)(968,473)(968,473)(969,473)(969,473)(969,473)(970,473)(970,473)(971,473)(971,473)(972,473)(972,473)(972,473)(973,473)(973,473)(974,473)(974,473)(975,473)
\psset{linewidth=0.35pt,linestyle=solid}\psline(975,473)(975,473)(975,473)(976,473)(976,473)(977,473)(977,473)(977,473)(978,473)(978,473)(979,473)(979,473)(980,473)(980,473)(980,473)(981,473)(981,473)(982,473)(982,473)(983,473)(983,473)(983,473)(984,473)(984,473)(985,473)(985,473)(985,473)(986,473)(986,473)(987,473)(987,473)(988,473)(988,473)(988,473)(989,473)(989,473)(990,473)(990,473)(991,473)(991,473)(991,473)(992,473)(992,473)(993,473)(993,473)(993,473)(994,473)(994,473)(995,473)(995,473)(996,473)
\psset{linewidth=0.35pt,linestyle=solid}\psline(996,473)(996,473)(996,473)(997,473)(997,473)(998,473)(998,473)(998,473)(999,473)(999,473)(1000,473)(1000,473)(1001,473)(1001,473)(1001,473)(1002,473)(1002,473)(1003,473)(1003,473)(1004,473)(1004,473)(1004,473)(1005,473)(1005,473)(1006,473)(1006,473)(1006,473)(1007,473)(1007,473)(1008,473)(1008,473)(1009,473)(1009,473)(1009,473)(1010,473)(1010,473)(1011,473)(1011,473)(1012,473)(1012,473)(1012,473)(1013,473)(1013,473)(1014,473)(1014,473)(1014,473)(1015,473)(1015,473)(1016,473)(1016,473)(1017,473)
\psset{linewidth=0.35pt,linestyle=solid}\psline(1017,473)(1017,473)
\end{pspicture}}
\caption{The function $\gamma$.}\label{gamma_figuur}
\end{figure}

Let $\gamma$ be the transformation on $[-2/\sqrt{1-q}, 2/\sqrt{1-q}]$
that changes the orientations of $[-2/\sqrt{1-q}, 0]$ and $[0,
2/\sqrt{1-q}]$ in such a way that the distribution $\nu_q$ is
preserved. For this $\gamma$ has to satisfy the differential equation:
\begin{equation}
\label{dv_voor_nu_q}
\nu_q^\prime(x) dx + \nu_q^\prime(\gamma(x)) d\gamma(x) 
=
	0,
\end{equation}
with $\gamma(-2/\sqrt{1-q}) = \gamma(2/\sqrt{1-q}) = 0$.  Indeed,
this leads to:
\[
\mathbb{P}(0 \leq T \leq x)
=
	\mathbb{P}(0 \leq \gamma(T) \leq x)
=
	\mathbb{P}(\gamma^{-1}(x) \leq T \leq 2/\sqrt{1-q}) 
\]
as can be seen by differentiating both sides with respect to $x$. 
Note that the function $\gamma$ is its own inverse.
Figure~\ref{gamma_figuur} shows a typical picture of the shape of the
function $\gamma$.

Let $\widehat{W}$ be the unitary operator on $L^2(\mathbb{R}, \nu_q)$ that
implements $\gamma$:
\[
(\widehat{W} f)(x)
=
	f(\gamma(x))\qquad\text{for $f \in L^2(\mathbb{R}, \nu_q)$.}
\]
This immediately implies that $\widehat{W}^2 = \unit$ since $\gamma
\circ \gamma = \mathrm{id}$, so $\widehat{W}$ is self-adjoint.
If we define $\widetilde{W} := U^*\widehat{W} U$ it follows that:
\[
\gamma(X) = \gamma(U^* T U) = U^* \gamma(T) U = U^* \widehat{W} T \widehat{W} U
= \widetilde{W} X \widetilde{W},
\]
so $\widetilde{W}$ is a unitary and self-adjoint operator on
$\mathcal{F}_q(\mathbb{C})$ that implements $\gamma$ on $X$.
Note that $\widetilde{W}\Omega=\Omega$ because 
$\widehat{W}1 = 1$. On the
canonical basis $(e_j)_{j \in \mathbb{N}}$ of $\mathcal{F}_q(\mathbb{C})$
the operator $\widetilde{W}$ can be written as:
\[
\widetilde{W} e_n = \sum_{k=0}^\infty w_{k n} e_k
\qquad\text{with $w_{0 0} = 1$.}
\]
Now let us choose $\mathcal{H} = \mathbb{C}^2$, $f_0 = (1, 0)$ and
$f_1 = (0, 1)$ and let us denote $a(f_0)$ by $a_0$ and $a(f_1)$ by
$a_1$. Recall from the introduction that $X_0 = a_0 + a_0^*$ and $X_1
= a_1 + a_1^*$. In this setting we need a unitary operator $W$ on
$\mathcal{F}_q(\mathbb{C}^2)$ that satisfies: $\gamma(X) = W X W$. To
this end we denote by $\mathcal{K} \subset
\mathcal{F}_q(\mathbb{C}^2)$ the kernel of the operator $a_0$ on
$\mathcal{F}_q(\mathbb{C}^2)$, then by constructing an isomorphism $V\colon
\mathcal{F}_q(\mathbb{C}) \otimes \mathcal{K}
\rightarrow \mathcal{F}_q(\mathbb{C}^2)$,
the operator $\widetilde{W}$ can be extended to $W = V (\widetilde{W}
\otimes \unit) V^*$.
\begin{proposition}
The space $\mathcal{F}_q(\mathbb{C}^2)$ is canonically isomorphic to
$\mathcal{F}_q(\mathbb{C})\otimes \mathcal{K}$.
\end{proposition}
\begin{proof}
From the commutation relation~(\ref{comm1}) we find that 
$a_0^n (a_0^*)^n = P_n (a_0^* a_0)$, where $P_n$ is a polynomial
of degree $n$ with constant coefficient $[n]_q!$. In fact, $P_n$
is given by 
\[
P_n(x) = \prod_{j=1}^{n}(q^{j} x + [j]_q).
\]

For $n \in \mathbb{N}$, let $\mathcal{K}_n$ denote the Hilbert subspace
$(a_0^*)^n\mathcal{K}$. Note that $\mathcal{K}_n$ is indeed closed, since
$(a_0^*)^n$ acts on $\mathcal{K}$ as a multiple of an isometry,
for every $\phi \in \mathcal{K}$,
\[
\| (a_0^*)^n\phi\|^2_q
=
        \langle \phi, P_n(a_0^* a_0)\phi \rangle_q
=
        [n]_q! \|\phi\|_q^2.
\]

Furthermore, $\mathcal{K}_n \perp \mathcal{K}_m$ for $n > m$, since
for $\phi, \psi \in \mathcal{K}$,
\begin{align*}
\langle (a_0^*)^n \phi, (a_0^*)^m \psi\rangle_q
&=
        \langle (a_0^*)^{n-m}\phi, a_0^m (a_0^*)^m \psi\rangle_q\\
&=
        \langle (a_0^*)^{n-m}\phi, P_m(a_0^* a_0)\psi\rangle_q\\
&=
        [m]_q! \langle \phi, (a_0)^{n-m}\psi \rangle_q = 0.
\end{align*}

Now suppose that some $\psi \in \mathcal{F}_q(\mathbb{C}^2)$ is orthogonal
to all the $\mathcal{K}_n$. We claim that for all $n \in \mathbb{N}$
\begin{equation}
\label{claim}
\psi \perp \ker a_0^n,
\end{equation}
from which it follows that $\psi = 0$, since
$(\mathbb{C}^2)^{\otimes n} = \mathcal{F}_q^{(n)}(\mathbb{C}^2)
\subset \ker a_0^{n+1}$.

We proceed to prove (\ref{claim}) by induction. For $n=1$ we already have
(\ref{claim}) since $\ker a_0 = \mathcal{K}$. Suppose that (\ref{claim})
holds for some $n$. Then $\psi \in \overline{\ran (a_0^*)^n}$, say
$\psi = \lim_{k \rightarrow \infty} (a_0^*)^n \phi_k$ with
$\phi_k \in \mathcal{F}_q(\mathbb{C}^2)$, $k \in \mathbb{N}$.
Take $\theta \in \ker a_0^{n+1}$ and define
$\xi := a_0^n \theta \in \mathcal{K}$, then
\begin{align*}
\langle \psi, \theta\rangle_q
&= \lim_{k\rightarrow\infty} \langle (a_0^*)^n\phi_k, \theta\rangle_q
=
        \lim_{k\rightarrow\infty} \langle \phi_k, \xi\rangle_q\\
&=
        \frac{1}{[n]_q!}\lim_{k\rightarrow\infty} \langle \phi_k,
        P_n(a_0^* a_0)\xi\rangle_q\\
&=
        \frac{1}{[n]_q!}\lim_{k\rightarrow\infty} \langle
        (a_0^*)^n\phi_k, (a_0^*)^n\xi\rangle_q\\
&=
        \frac{1}{[n]_q!} \langle \psi, (a_0^*)^n\xi\rangle_q = 0,
\end{align*}
because $(a_0^*)^n\xi \in \mathcal{K}_n \perp \psi$. The claim,
(\ref{claim}), follows by induction.

We define an operator $V\colon
\mathcal{F}_q(\mathbb{C})\otimes\mathcal{K}
\rightarrow \mathcal{F}_q(\mathbb{C}^2)$ by:
\[
V(e_n \otimes \phi)
:=
        (a_0^*)^n \phi.
\]
The operator $V$ is an isomorphism since its range is dense by the
above, and, for all $\phi, \xi \in \mathcal{K}$,
\begin{align*}
\langle V(e_n \otimes \phi),
V(e_m \otimes \xi) \rangle_q
&=
        \langle (a_0^*)^n\phi, (a_0^*)^m\xi\rangle_q\\
&=
        \delta_{n, m} \langle \phi, P_n(a_0^* a_0)\xi\rangle_q\\
&=
        \delta_{n, m} \lbrack n\rbrack_q! \langle \phi, \xi\rangle_q\\
&=
        \langle e_n \otimes \phi,
        e_m \otimes \xi\rangle_q.
\end{align*}
\end{proof}
\begin{lemma}
The operator $W$ has the following properties:
\begin{enumerate}
\item{$W$ is unitary and self adjoint,}
\label{lemma1:1}
\item{$\gamma(X_0) = W X_0 W$,}
\label{lemma1:2}
\item{$W\phi = \phi$ for all $\phi \in \mathcal{K}$, in particular
      $W\Omega = \Omega$,}
\label{lemma1:3}
\item{$W(X_0 \phi) = \sum_{k = 1}^\infty w_{k1} (a_0^*)^k \phi$ for
      all $\phi \in \mathcal{K}$.}
\label{lemma1:4}
\end{enumerate}
\end{lemma}
\begin{proof}
Property \ref{lemma1:1} is clear from the definition of $W$ since
$\widetilde{W}$ is unitary and self adjoint.

To prove property \ref{lemma1:2}, note that, for $\phi \in \mathcal{K}$
and $n \in \mathbb{N}$, $V (a^* \otimes \unit) (e_n \otimes \phi)
= V (e_{n+1} \otimes \phi) = a_0^* V (e_n \otimes \phi)$, so
$V (a^* \otimes \unit) V^* = a_0^*$ and $V (X \otimes \unit) V^* = X_0$.
It follows that
\begin{align*}
W X_0 W 
&=
        W V (X \otimes \unit) V^* W\\
&=
        V (\widetilde{W} \otimes \unit) (X \otimes \unit)
        (\widetilde{W} \otimes \unit) V^*\\
&=
        V (\gamma(X) \otimes \unit) V^* \\
&=
        V \gamma(X \otimes \unit) V^* = \gamma (X_0).
\end{align*}

Property \ref{lemma1:3} is immediate from definitions:
\[
W \phi
=
        V (\widetilde{W} \otimes \unit) (e_0 \otimes \phi)
=
        V (e_0 \otimes \phi) = \phi,
\]
for all $\phi \in \mathcal{K}$.

The proof of property \ref{lemma1:4} is also immediate from definitions:
\begin{align*}
W (X_0 \phi)
&=
        W(a_0^* \phi)
=
        V (\widetilde{W} \otimes \unit) (e_1 \otimes \phi)
=
        V (\widetilde{W} e_1 \otimes \phi)\\
&=
        \sum_{k=1}^\infty w_{k 1} V (e_k \otimes \phi)
=
        \sum_{k=1}^\infty w_{k 1} (a_0^*)^k\phi.
\end{align*}
\end{proof}

We now turn to the proof of theorem~\ref{hoofdstelling}.
\begin{proof}
First we calculate the fourth moment of $X_0+X_1$. Since $X(f_0 + f_1)$
is $q$-Gaussian with variance $2$ we have:
\begin{align*}
\langle \Omega, (X_0+X_1)^4\Omega\rangle_q
&=
        (\sqrt{2})^4 \langle \Omega, X_0^4\Omega\rangle_q
=
        4 \| X_0^2\Omega \|_q^2 \\
&=
        4 (\|\Omega\|_q^2 + \|f^{\otimes 2}\|_q^2)
=
        4 (1 + [2]_q)\\
&=
        8 + 4 q,
\end{align*}
a linear interpolation between $8$ and $12$ for $q$ varying between
$0$ and $1$.
We now turn to the calculation of the fourth moment of $\gamma(X_0)+X_1$:
\[
\langle \Omega, (\gamma(X_0)+X_1)^4\Omega\rangle_q
=
        \| (\gamma(X_0)+X_1)^2\Omega \|_q^2.
\]
For this we need the following:
\begin{align*}
\gamma(X_0)^2 \Omega
&=
        W X_0^2 W \Omega = \Omega + W f_0^{\otimes 2}\\
X_1^2 \Omega
&=
        \Omega + f_1^{\otimes 2}\\
\gamma(X_0) X_1 \Omega
&=
        W X_0 W X_1 \Omega = W X_0 X_1 \Omega\\
&=
        \sum_{k=1}^\infty w_{k 1} (a_0^*)^k X_1 \Omega =
        \sum_{k=1}^\infty w_{k 1} f_0^{\otimes k} \otimes f_1\\
X_1 \gamma(X_0) \Omega
&=
        X_1 W X_0 \Omega =
        \sum_{k=1}^\infty w_{k 1} f_1 \otimes f_0^{\otimes k},
\end{align*}
from which it is easy to deduce that:
\begin{equation}
\label{bewijs3:1}
\| (\gamma(X_0) + X_1)^2 \Omega\|_q^2
=
        \| (\gamma(X_0)^2 + X_1^2) \Omega\|_q^2 +
        \| (\gamma(X_0) X_1 + X_1 \gamma(X_0)) \Omega \|_q^2.
\end{equation}
The first term on the right hand side of (\ref{bewijs3:1}) is found
to be:
\[
\| (\gamma(X_0)^2 + X_1^2) \Omega\|_q^2
=
        4 \| \Omega \|_q^2 + \| f_0^{\otimes 2} \|_q^2 +
        \| f_1^{\otimes 2} \|_q^2 = 4 + 2 [2]_q = 6 + 2 q.
\]
\begin{figure}
\centerline{
\psset{unit=0.240900pt}
\psset{arrowsize=3pt 3.2 1.4 .3}
\begin{pspicture}(0,0)(1081,900)
\psset{linewidth=0.3pt,linestyle=dotted,dotsep=1pt}\psline(176,68)(176,877)
\psset{linewidth=0.35pt,linestyle=solid}\psline(176,68)(196,68)
\psset{linewidth=0.35pt,linestyle=solid}\psline(1017,68)(997,68)
\rput[r](154,68){$8$}
\psset{linewidth=0.35pt,linestyle=solid}\psline(176,270)(196,270)
\psset{linewidth=0.35pt,linestyle=solid}\psline(1017,270)(997,270)
\rput[r](154,270){$9$}
\psset{linewidth=0.35pt,linestyle=solid}\psline(176,473)(196,473)
\psset{linewidth=0.35pt,linestyle=solid}\psline(1017,473)(997,473)
\rput[r](154,473){$10$}
\psset{linewidth=0.35pt,linestyle=solid}\psline(176,675)(196,675)
\psset{linewidth=0.35pt,linestyle=solid}\psline(1017,675)(997,675)
\rput[r](154,675){$11$}
\psset{linewidth=0.35pt,linestyle=solid}\psline(176,877)(196,877)
\psset{linewidth=0.35pt,linestyle=solid}\psline(1017,877)(997,877)
\rput[r](154,877){$12$}
\psset{linewidth=0.35pt,linestyle=solid}\psline(176,68)(176,88)
\psset{linewidth=0.35pt,linestyle=solid}\psline(176,877)(176,857)
\rput(176,23){0}
\psset{linewidth=0.35pt,linestyle=solid}\psline(344,68)(344,88)
\psset{linewidth=0.35pt,linestyle=solid}\psline(344,877)(344,857)
\rput(344,23){0.2}
\psset{linewidth=0.35pt,linestyle=solid}\psline(512,68)(512,88)
\psset{linewidth=0.35pt,linestyle=solid}\psline(512,877)(512,857)
\rput(512,23){0.4}
\psset{linewidth=0.35pt,linestyle=solid}\psline(681,68)(681,88)
\psset{linewidth=0.35pt,linestyle=solid}\psline(681,877)(681,857)
\rput(681,23){0.6}
\psset{linewidth=0.35pt,linestyle=solid}\psline(849,68)(849,88)
\psset{linewidth=0.35pt,linestyle=solid}\psline(849,877)(849,857)
\rput(849,23){0.8}
\psset{linewidth=0.35pt,linestyle=solid}\psline(1017,68)(1017,88)
\psset{linewidth=0.35pt,linestyle=solid}\psline(1017,877)(1017,857)
\rput(1017,23){1}
\psset{linewidth=0.35pt,linestyle=solid}\psline(176,68)(1017,68)(1017,877)(176,877)(176,68)
\psset{linewidth=0.35pt,linestyle=solid}\psline(176,68)(176,68)(597,473)(1017,877)
\psset{linewidth=0.35pt,linestyle=solid}\psline(176,68)(176,68)(176,68)(177,68)(177,69)(178,69)(178,69)(179,69)(179,70)(179,70)(180,70)(180,70)(181,71)(181,71)(181,71)(182,71)(182,72)(183,72)(183,72)(184,72)(184,73)(184,73)(185,73)(185,73)(186,74)(186,74)(187,74)(187,74)(187,75)(188,75)(188,75)(189,75)(189,75)(189,76)(190,76)(190,76)(191,76)(191,77)(192,77)(192,77)(192,77)(193,78)(193,78)(194,78)(194,78)(195,79)(195,79)(195,79)(196,79)(196,80)(197,80)
\psset{linewidth=0.35pt,linestyle=solid}\psline(197,80)(197,80)(197,80)(198,81)(198,81)(199,81)(199,81)(200,81)(200,82)(200,82)(201,82)(201,82)(202,83)(202,83)(202,83)(203,83)(203,84)(204,84)(204,84)(205,84)(205,85)(205,85)(206,85)(206,85)(207,86)(207,86)(208,86)(208,86)(208,87)(209,87)(209,87)(210,87)(210,87)(210,88)(211,88)(211,88)(212,88)(212,89)(213,89)(213,89)(213,89)(214,90)(214,90)(215,90)(215,90)(216,91)(216,91)(216,91)(217,91)(217,92)(218,92)
\psset{linewidth=0.35pt,linestyle=solid}\psline(218,92)(218,92)(218,92)(219,93)(219,93)(220,93)(220,93)(221,93)(221,94)(221,94)(222,94)(222,94)(223,95)(223,95)(224,95)(224,95)(224,96)(225,96)(225,96)(226,96)(226,97)(226,97)(227,97)(227,97)(228,98)(228,98)(229,98)(229,98)(229,98)(230,99)(230,99)(231,99)(231,99)(232,100)(232,100)(232,100)(233,100)(233,101)(234,101)(234,101)(234,101)(235,102)(235,102)(236,102)(236,102)(237,102)(237,103)(237,103)(238,103)(238,103)(239,104)
\psset{linewidth=0.35pt,linestyle=solid}\psline(239,104)(239,104)(239,104)(240,104)(240,105)(241,105)(241,105)(242,105)(242,106)(242,106)(243,106)(243,106)(244,107)(244,107)(245,107)(245,107)(245,107)(246,108)(246,108)(247,108)(247,108)(247,109)(248,109)(248,109)(249,109)(249,110)(250,110)(250,110)(250,110)(251,111)(251,111)(252,111)(252,111)(253,111)(253,112)(253,112)(254,112)(254,112)(255,113)(255,113)(255,113)(256,113)(256,114)(257,114)(257,114)(258,114)(258,115)(258,115)(259,115)(259,115)(260,116)
\psset{linewidth=0.35pt,linestyle=solid}\psline(260,116)(260,116)(261,116)(261,116)(261,116)(262,117)(262,117)(263,117)(263,117)(263,118)(264,118)(264,118)(265,118)(265,119)(266,119)(266,119)(266,119)(267,120)(267,120)(268,120)(268,120)(269,120)(269,121)(269,121)(270,121)(270,121)(271,122)(271,122)(271,122)(272,122)(272,123)(273,123)(273,123)(274,123)(274,124)(274,124)(275,124)(275,124)(276,124)(276,125)(276,125)(277,125)(277,125)(278,126)(278,126)(279,126)(279,126)(279,127)(280,127)(280,127)(281,127)
\psset{linewidth=0.35pt,linestyle=solid}\psline(281,127)(281,128)(282,128)(282,128)(282,128)(283,128)(283,129)(284,129)(284,129)(284,129)(285,130)(285,130)(286,130)(286,130)(287,131)(287,131)(287,131)(288,131)(288,132)(289,132)(289,132)(290,132)(290,132)(290,133)(291,133)(291,133)(292,133)(292,134)(292,134)(293,134)(293,134)(294,135)(294,135)(295,135)(295,135)(295,136)(296,136)(296,136)(297,136)(297,136)(298,137)(298,137)(298,137)(299,137)(299,138)(300,138)(300,138)(300,138)(301,139)(301,139)(302,139)
\psset{linewidth=0.35pt,linestyle=solid}\psline(302,139)(302,139)(303,140)(303,140)(303,140)(304,140)(304,140)(305,141)(305,141)(306,141)(306,141)(306,142)(307,142)(307,142)(308,142)(308,143)(308,143)(309,143)(309,143)(310,144)(310,144)(311,144)(311,144)(311,144)(312,145)(312,145)(313,145)(313,145)(314,146)(314,146)(314,146)(315,146)(315,147)(316,147)(316,147)(316,147)(317,148)(317,148)(318,148)(318,148)(319,148)(319,149)(319,149)(320,149)(320,149)(321,150)(321,150)(321,150)(322,150)(322,151)(323,151)
\psset{linewidth=0.35pt,linestyle=solid}\psline(323,151)(323,151)(324,151)(324,152)(324,152)(325,152)(325,152)(326,152)(326,153)(327,153)(327,153)(327,153)(328,154)(328,154)(329,154)(329,154)(329,155)(330,155)(330,155)(331,155)(331,155)(332,156)(332,156)(332,156)(333,156)(333,157)(334,157)(334,157)(335,157)(335,158)(335,158)(336,158)(336,158)(337,159)(337,159)(337,159)(338,159)(338,159)(339,160)(339,160)(340,160)(340,160)(340,161)(341,161)(341,161)(342,161)(342,162)(343,162)(343,162)(343,162)(344,163)
\psset{linewidth=0.35pt,linestyle=solid}\psline(344,163)(344,163)(345,163)(345,163)(345,163)(346,164)(346,164)(347,164)(347,164)(348,165)(348,165)(348,165)(349,165)(349,166)(350,166)(350,166)(351,166)(351,167)(351,167)(352,167)(352,167)(353,167)(353,168)(353,168)(354,168)(354,168)(355,169)(355,169)(356,169)(356,169)(356,170)(357,170)(357,170)(358,170)(358,171)(358,171)(359,171)(359,171)(360,171)(360,172)(361,172)(361,172)(361,172)(362,173)(362,173)(363,173)(363,173)(364,174)(364,174)(364,174)(365,174)
\psset{linewidth=0.35pt,linestyle=solid}\psline(365,174)(365,175)(366,175)(366,175)(366,175)(367,175)(367,176)(368,176)(368,176)(369,176)(369,177)(369,177)(370,177)(370,177)(371,178)(371,178)(372,178)(372,178)(372,178)(373,179)(373,179)(374,179)(374,179)(374,180)(375,180)(375,180)(376,180)(376,181)(377,181)(377,181)(377,181)(378,182)(378,182)(379,182)(379,182)(380,182)(380,183)(380,183)(381,183)(381,183)(382,184)(382,184)(382,184)(383,184)(383,185)(384,185)(384,185)(385,185)(385,186)(385,186)(386,186)
\psset{linewidth=0.35pt,linestyle=solid}\psline(386,186)(386,186)(387,186)(387,187)(388,187)(388,187)(388,187)(389,188)(389,188)(390,188)(390,188)(390,189)(391,189)(391,189)(392,189)(392,190)(393,190)(393,190)(393,190)(394,190)(394,191)(395,191)(395,191)(396,191)(396,192)(396,192)(397,192)(397,192)(398,193)(398,193)(398,193)(399,193)(399,194)(400,194)(400,194)(401,194)(401,194)(401,195)(402,195)(402,195)(403,195)(403,196)(403,196)(404,196)(404,196)(405,197)(405,197)(406,197)(406,197)(406,198)(407,198)
\psset{linewidth=0.35pt,linestyle=solid}\psline(407,198)(407,198)(408,198)(408,198)(409,199)(409,199)(409,199)(410,199)(410,200)(411,200)(411,200)(411,200)(412,201)(412,201)(413,201)(413,201)(414,202)(414,202)(414,202)(415,202)(415,203)(416,203)(416,203)(417,203)(417,203)(417,204)(418,204)(418,204)(419,204)(419,205)(419,205)(420,205)(420,205)(421,206)(421,206)(422,206)(422,206)(422,207)(423,207)(423,207)(424,207)(424,207)(425,208)(425,208)(425,208)(426,208)(426,209)(427,209)(427,209)(427,209)(428,210)
\psset{linewidth=0.35pt,linestyle=solid}\psline(428,210)(428,210)(429,210)(429,210)(430,211)(430,211)(430,211)(431,211)(431,211)(432,212)(432,212)(433,212)(433,212)(433,213)(434,213)(434,213)(435,213)(435,214)(435,214)(436,214)(436,214)(437,215)(437,215)(438,215)(438,215)(438,216)(439,216)(439,216)(440,216)(440,216)(440,217)(441,217)(441,217)(442,217)(442,218)(443,218)(443,218)(443,218)(444,219)(444,219)(445,219)(445,219)(446,220)(446,220)(446,220)(447,220)(447,220)(448,221)(448,221)(448,221)(449,221)
\psset{linewidth=0.35pt,linestyle=solid}\psline(449,221)(449,222)(450,222)(450,222)(451,222)(451,223)(451,223)(452,223)(452,223)(453,224)(453,224)(454,224)(454,224)(454,225)(455,225)(455,225)(456,225)(456,225)(456,226)(457,226)(457,226)(458,226)(458,227)(459,227)(459,227)(459,227)(460,228)(460,228)(461,228)(461,228)(462,229)(462,229)(462,229)(463,229)(463,230)(464,230)(464,230)(464,230)(465,230)(465,231)(466,231)(466,231)(467,231)(467,232)(467,232)(468,232)(468,232)(469,233)(469,233)(470,233)(470,233)
\psset{linewidth=0.35pt,linestyle=solid}\psline(470,233)(470,234)(471,234)(471,234)(472,234)(472,235)(472,235)(473,235)(473,235)(474,235)(474,236)(475,236)(475,236)(475,236)(476,237)(476,237)(477,237)(477,237)(477,238)(478,238)(478,238)(479,238)(479,239)(480,239)(480,239)(480,239)(481,240)(481,240)(482,240)(482,240)(483,240)(483,241)(483,241)(484,241)(484,241)(485,242)(485,242)(485,242)(486,242)(486,243)(487,243)(487,243)(488,243)(488,244)(488,244)(489,244)(489,244)(490,245)(490,245)(491,245)(491,245)
\psset{linewidth=0.35pt,linestyle=solid}\psline(491,245)(491,246)(492,246)(492,246)(493,246)(493,246)(493,247)(494,247)(494,247)(495,247)(495,248)(496,248)(496,248)(496,248)(497,249)(497,249)(498,249)(498,249)(499,250)(499,250)(499,250)(500,250)(500,251)(501,251)(501,251)(501,251)(502,252)(502,252)(503,252)(503,252)(504,252)(504,253)(504,253)(505,253)(505,253)(506,254)(506,254)(507,254)(507,254)(507,255)(508,255)(508,255)(509,255)(509,256)(509,256)(510,256)(510,256)(511,257)(511,257)(512,257)(512,257)
\psset{linewidth=0.35pt,linestyle=solid}\psline(512,257)(512,258)(513,258)(513,258)(514,258)(514,258)(515,259)(515,259)(515,259)(516,259)(516,260)(517,260)(517,260)(517,260)(518,261)(518,261)(519,261)(519,261)(520,262)(520,262)(520,262)(521,262)(521,263)(522,263)(522,263)(522,263)(523,264)(523,264)(524,264)(524,264)(525,265)(525,265)(525,265)(526,265)(526,266)(527,266)(527,266)(528,266)(528,266)(528,267)(529,267)(529,267)(530,267)(530,268)(530,268)(531,268)(531,268)(532,269)(532,269)(533,269)(533,269)
\psset{linewidth=0.35pt,linestyle=solid}\psline(533,269)(533,270)(534,270)(534,270)(535,270)(535,271)(536,271)(536,271)(536,271)(537,272)(537,272)(538,272)(538,272)(538,273)(539,273)(539,273)(540,273)(540,274)(541,274)(541,274)(541,274)(542,274)(542,275)(543,275)(543,275)(544,275)(544,276)(544,276)(545,276)(545,276)(546,277)(546,277)(546,277)(547,277)(547,278)(548,278)(548,278)(549,278)(549,279)(549,279)(550,279)(550,279)(551,280)(551,280)(552,280)(552,280)(552,281)(553,281)(553,281)(554,281)(554,282)
\psset{linewidth=0.35pt,linestyle=solid}\psline(554,282)(554,282)(555,282)(555,282)(556,283)(556,283)(557,283)(557,283)(557,284)(558,284)(558,284)(559,284)(559,285)(559,285)(560,285)(560,285)(561,285)(561,286)(562,286)(562,286)(562,286)(563,287)(563,287)(564,287)(564,287)(565,288)(565,288)(565,288)(566,288)(566,289)(567,289)(567,289)(567,289)(568,290)(568,290)(569,290)(569,290)(570,291)(570,291)(570,291)(571,291)(571,292)(572,292)(572,292)(573,292)(573,293)(573,293)(574,293)(574,293)(575,294)(575,294)
\psset{linewidth=0.35pt,linestyle=solid}\psline(575,294)(575,294)(576,294)(576,295)(577,295)(577,295)(578,295)(578,296)(578,296)(579,296)(579,296)(580,297)(580,297)(581,297)(581,297)(581,298)(582,298)(582,298)(583,298)(583,299)(583,299)(584,299)(584,299)(585,300)(585,300)(586,300)(586,300)(586,301)(587,301)(587,301)(588,301)(588,302)(589,302)(589,302)(589,302)(590,303)(590,303)(591,303)(591,303)(591,304)(592,304)(592,304)(593,304)(593,305)(594,305)(594,305)(594,305)(595,306)(595,306)(596,306)(596,306)
\psset{linewidth=0.35pt,linestyle=solid}\psline(596,306)(597,307)(597,307)(597,307)(598,307)(598,308)(599,308)(599,308)(599,308)(600,308)(600,309)(601,309)(601,309)(602,309)(602,310)(602,310)(603,310)(603,311)(604,311)(604,311)(604,311)(605,312)(605,312)(606,312)(606,312)(607,313)(607,313)(607,313)(608,313)(608,314)(609,314)(609,314)(610,314)(610,315)(610,315)(611,315)(611,315)(612,316)(612,316)(612,316)(613,316)(613,317)(614,317)(614,317)(615,317)(615,318)(615,318)(616,318)(616,318)(617,319)(617,319)
\psset{linewidth=0.35pt,linestyle=solid}\psline(617,319)(618,319)(618,319)(618,320)(619,320)(619,320)(620,320)(620,321)(620,321)(621,321)(621,321)(622,322)(622,322)(623,322)(623,322)(623,323)(624,323)(624,323)(625,323)(625,324)(626,324)(626,324)(626,324)(627,325)(627,325)(628,325)(628,325)(628,326)(629,326)(629,326)(630,326)(630,327)(631,327)(631,327)(631,327)(632,328)(632,328)(633,328)(633,328)(634,329)(634,329)(634,329)(635,329)(635,330)(636,330)(636,330)(636,330)(637,331)(637,331)(638,331)(638,331)
\psset{linewidth=0.35pt,linestyle=solid}\psline(638,331)(639,332)(639,332)(639,332)(640,333)(640,333)(641,333)(641,333)(641,334)(642,334)(642,334)(643,334)(643,335)(644,335)(644,335)(644,335)(645,336)(645,336)(646,336)(646,336)(647,337)(647,337)(647,337)(648,337)(648,338)(649,338)(649,338)(649,338)(650,339)(650,339)(651,339)(651,339)(652,340)(652,340)(652,340)(653,340)(653,341)(654,341)(654,341)(655,342)(655,342)(655,342)(656,342)(656,343)(657,343)(657,343)(657,343)(658,344)(658,344)(659,344)(659,344)
\psset{linewidth=0.35pt,linestyle=solid}\psline(659,344)(660,345)(660,345)(660,345)(661,345)(661,346)(662,346)(662,346)(663,346)(663,347)(663,347)(664,347)(664,347)(665,348)(665,348)(665,348)(666,349)(666,349)(667,349)(667,349)(668,350)(668,350)(668,350)(669,350)(669,351)(670,351)(670,351)(671,351)(671,352)(671,352)(672,352)(672,352)(673,353)(673,353)(673,353)(674,353)(674,354)(675,354)(675,354)(676,355)(676,355)(676,355)(677,355)(677,356)(678,356)(678,356)(678,356)(679,357)(679,357)(680,357)(680,357)
\psset{linewidth=0.35pt,linestyle=solid}\psline(680,357)(681,358)(681,358)(681,358)(682,358)(682,359)(683,359)(683,359)(684,360)(684,360)(684,360)(685,360)(685,361)(686,361)(686,361)(686,361)(687,362)(687,362)(688,362)(688,362)(689,363)(689,363)(689,363)(690,364)(690,364)(691,364)(691,364)(692,365)(692,365)(692,365)(693,365)(693,366)(694,366)(694,366)(694,366)(695,367)(695,367)(696,367)(696,368)(697,368)(697,368)(697,368)(698,369)(698,369)(699,369)(699,369)(700,370)(700,370)(700,370)(701,370)(701,371)
\psset{linewidth=0.35pt,linestyle=solid}\psline(701,371)(702,371)(702,371)(702,372)(703,372)(703,372)(704,372)(704,373)(705,373)(705,373)(705,373)(706,374)(706,374)(707,374)(707,374)(708,375)(708,375)(708,375)(709,376)(709,376)(710,376)(710,376)(710,377)(711,377)(711,377)(712,377)(712,378)(713,378)(713,378)(713,379)(714,379)(714,379)(715,379)(715,380)(716,380)(716,380)(716,380)(717,381)(717,381)(718,381)(718,382)(718,382)(719,382)(719,382)(720,383)(720,383)(721,383)(721,383)(721,384)(722,384)(722,384)
\psset{linewidth=0.35pt,linestyle=solid}\psline(722,384)(723,385)(723,385)(723,385)(724,385)(724,386)(725,386)(725,386)(726,386)(726,387)(726,387)(727,387)(727,388)(728,388)(728,388)(729,388)(729,389)(729,389)(730,389)(730,390)(731,390)(731,390)(731,390)(732,391)(732,391)(733,391)(733,391)(734,392)(734,392)(734,392)(735,393)(735,393)(736,393)(736,393)(737,394)(737,394)(737,394)(738,395)(738,395)(739,395)(739,395)(739,396)(740,396)(740,396)(741,396)(741,397)(742,397)(742,397)(742,398)(743,398)(743,398)
\psset{linewidth=0.35pt,linestyle=solid}\psline(743,398)(744,398)(744,399)(745,399)(745,399)(745,400)(746,400)(746,400)(747,400)(747,401)(747,401)(748,401)(748,402)(749,402)(749,402)(750,402)(750,403)(750,403)(751,403)(751,404)(752,404)(752,404)(753,404)(753,405)(753,405)(754,405)(754,406)(755,406)(755,406)(755,406)(756,407)(756,407)(757,407)(757,408)(758,408)(758,408)(758,408)(759,409)(759,409)(760,409)(760,410)(760,410)(761,410)(761,410)(762,411)(762,411)(763,411)(763,412)(763,412)(764,412)(764,412)
\psset{linewidth=0.35pt,linestyle=solid}\psline(764,412)(765,413)(765,413)(766,413)(766,414)(766,414)(767,414)(767,414)(768,415)(768,415)(768,415)(769,416)(769,416)(770,416)(770,416)(771,417)(771,417)(771,417)(772,418)(772,418)(773,418)(773,419)(774,419)(774,419)(774,419)(775,420)(775,420)(776,420)(776,421)(776,421)(777,421)(777,421)(778,422)(778,422)(779,422)(779,423)(779,423)(780,423)(780,424)(781,424)(781,424)(782,424)(782,425)(782,425)(783,425)(783,426)(784,426)(784,426)(784,426)(785,427)(785,427)
\psset{linewidth=0.35pt,linestyle=solid}\psline(785,427)(786,427)(786,428)(787,428)(787,428)(787,429)(788,429)(788,429)(789,429)(789,430)(790,430)(790,430)(790,431)(791,431)(791,431)(792,432)(792,432)(792,432)(793,433)(793,433)(794,433)(794,433)(795,434)(795,434)(795,434)(796,435)(796,435)(797,435)(797,436)(797,436)(798,436)(798,436)(799,437)(799,437)(800,437)(800,438)(800,438)(801,438)(801,439)(802,439)(802,439)(803,440)(803,440)(803,440)(804,440)(804,441)(805,441)(805,441)(805,442)(806,442)(806,442)
\psset{linewidth=0.35pt,linestyle=solid}\psline(806,442)(807,443)(807,443)(808,443)(808,444)(808,444)(809,444)(809,444)(810,445)(810,445)(811,445)(811,446)(811,446)(812,446)(812,447)(813,447)(813,447)(813,448)(814,448)(814,448)(815,449)(815,449)(816,449)(816,450)(816,450)(817,450)(817,450)(818,451)(818,451)(819,451)(819,452)(819,452)(820,452)(820,453)(821,453)(821,453)(821,454)(822,454)(822,454)(823,455)(823,455)(824,455)(824,456)(824,456)(825,456)(825,457)(826,457)(826,457)(827,458)(827,458)(827,458)
\psset{linewidth=0.35pt,linestyle=solid}\psline(827,458)(828,459)(828,459)(829,459)(829,459)(829,460)(830,460)(830,460)(831,461)(831,461)(832,461)(832,462)(832,462)(833,462)(833,463)(834,463)(834,463)(835,464)(835,464)(835,464)(836,465)(836,465)(837,465)(837,466)(837,466)(838,466)(838,467)(839,467)(839,467)(840,468)(840,468)(840,468)(841,469)(841,469)(842,469)(842,470)(842,470)(843,470)(843,471)(844,471)(844,471)(845,472)(845,472)(845,472)(846,473)(846,473)(847,474)(847,474)(848,474)(848,475)(848,475)
\psset{linewidth=0.35pt,linestyle=solid}\psline(848,475)(849,475)(849,476)(850,476)(850,476)(850,477)(851,477)(851,477)(852,478)(852,478)(853,478)(853,479)(853,479)(854,479)(854,480)(855,480)(855,480)(856,481)(856,481)(856,481)(857,482)(857,482)(858,483)(858,483)(858,483)(859,484)(859,484)(860,484)(860,485)(861,485)(861,485)(861,486)(862,486)(862,486)(863,487)(863,487)(864,488)(864,488)(864,488)(865,489)(865,489)(866,489)(866,490)(866,490)(867,490)(867,491)(868,491)(868,492)(869,492)(869,492)(869,493)
\psset{linewidth=0.35pt,linestyle=solid}\psline(869,493)(870,493)(870,493)(871,494)(871,494)(872,494)(872,495)(872,495)(873,496)(873,496)(874,496)(874,497)(874,497)(875,497)(875,498)(876,498)(876,499)(877,499)(877,499)(877,500)(878,500)(878,500)(879,501)(879,501)(879,502)(880,502)(880,502)(881,503)(881,503)(882,504)(882,504)(882,504)(883,505)(883,505)(884,505)(884,506)(885,506)(885,507)(885,507)(886,507)(886,508)(887,508)(887,509)(887,509)(888,509)(888,510)(889,510)(889,511)(890,511)(890,511)(890,512)
\psset{linewidth=0.35pt,linestyle=solid}\psline(890,512)(891,512)(891,513)(892,513)(892,513)(893,514)(893,514)(893,515)(894,515)(894,515)(895,516)(895,516)(895,517)(896,517)(896,517)(897,518)(897,518)(898,519)(898,519)(898,520)(899,520)(899,520)(900,521)(900,521)(901,522)(901,522)(901,522)(902,523)(902,523)(903,524)(903,524)(903,525)(904,525)(904,525)(905,526)(905,526)(906,527)(906,527)(906,528)(907,528)(907,528)(908,529)(908,529)(909,530)(909,530)(909,531)(910,531)(910,532)(911,532)(911,532)(911,533)
\psset{linewidth=0.35pt,linestyle=solid}\psline(911,533)(912,533)(912,534)(913,534)(913,535)(914,535)(914,536)(914,536)(915,536)(915,537)(916,537)(916,538)(917,538)(917,539)(917,539)(918,540)(918,540)(919,541)(919,541)(919,542)(920,542)(920,542)(921,543)(921,543)(922,544)(922,544)(922,545)(923,545)(923,546)(924,546)(924,547)(924,547)(925,548)(925,548)(926,549)(926,549)(927,550)(927,550)(927,551)(928,551)(928,552)(929,552)(929,553)(930,553)(930,554)(930,554)(931,555)(931,555)(932,556)(932,556)(932,557)
\psset{linewidth=0.35pt,linestyle=solid}\psline(932,557)(933,557)(933,558)(934,558)(934,559)(935,559)(935,560)(935,560)(936,561)(936,562)(937,562)(937,563)(938,563)(938,564)(938,564)(939,565)(939,565)(940,566)(940,566)(940,567)(941,568)(941,568)(942,569)(942,569)(943,570)(943,570)(943,571)(944,572)(944,572)(945,573)(945,573)(946,574)(946,574)(946,575)(947,576)(947,576)(948,577)(948,577)(948,578)(949,579)(949,579)(950,580)(950,580)(951,581)(951,582)(951,582)(952,583)(952,584)(953,584)(953,585)(954,586)
\psset{linewidth=0.35pt,linestyle=solid}\psline(954,586)(954,586)(954,587)(955,588)(955,589)(956,589)(956,590)(956,591)(957,592)(957,592)(958,593)(958,594)(959,595)(959,595)(959,596)(960,597)(960,598)(961,598)(961,599)(961,600)(962,601)(962,602)(963,602)(963,603)(964,604)(964,605)(964,606)(965,606)(965,607)(966,608)(966,609)(967,610)(967,610)(967,611)(968,612)(968,613)(969,614)(969,614)(969,615)(970,616)(970,617)(971,618)(971,618)(972,619)(972,620)(972,621)(973,622)(973,622)(974,623)(974,624)(975,625)
\psset{linewidth=0.35pt,linestyle=solid}\psline(975,625)(975,626)(975,627)(976,627)(976,628)(977,629)(977,630)(977,631)(978,632)(978,633)(979,633)(979,634)(980,635)(980,636)(980,637)(981,638)(981,639)(982,640)(982,641)(983,642)(983,643)(983,644)(984,646)(984,647)(985,648)(985,649)(985,650)(986,651)(986,653)(987,654)(987,655)(988,657)(988,658)(988,659)(989,661)(989,662)(990,664)(990,666)(991,667)(991,669)(991,670)(992,672)(992,674)(993,676)(993,678)(993,679)(994,681)(994,683)(995,685)(995,687)(996,689)
\psset{linewidth=0.35pt,linestyle=solid}\psline(996,689)(996,692)(996,694)(997,696)(997,698)(998,701)(998,703)(998,705)(999,708)(999,710)(1000,713)(1000,715)(1001,718)(1001,721)(1001,723)(1002,726)(1002,729)(1003,732)(1003,734)(1004,737)(1004,740)(1004,743)(1005,746)(1005,749)(1006,753)(1006,753)(1006,759)(1007,761)(1007,766)(1008,768)(1008,770)(1009,776)(1009,781)(1009,783)(1010,788)(1010,791)(1011,795)(1011,800)(1012,804)(1012,808)(1012,814)(1013,818)(1013,823)(1014,828)(1014,833)(1014,838)(1015,845)(1015,850)(1016,856)(1016,862)(1017,870)
\psset{linewidth=0.35pt,linestyle=solid}\psline(1017,870)(1017,877)
\end{pspicture}}
\caption{The fourth moment of $X+Y$ and $\gamma(X)+Y$ for $q \in (0,1)$.}
\label{mom_figuur}
\end{figure}
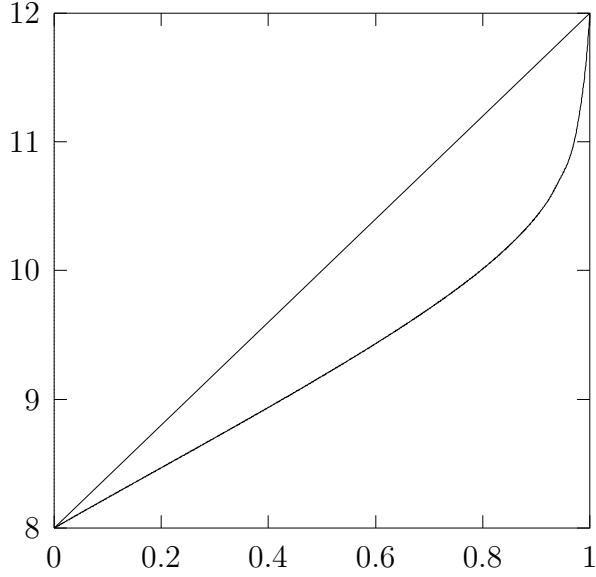
The second term on the right hand side of (\ref{bewijs3:1}) yields:
\begin{align*}
\| (\gamma(X_0) X_1 + X_1 \gamma(X_0)) \Omega \|_q^2
&=
        \sum_{k=1}^\infty w_{k 1}^2 \|
        (f_0^{\otimes k} \otimes f_1 + f_1 \otimes f_0^{\otimes k} ) \|_q^2\\
&=
        \sum_{k=1}^\infty w_{k 1}^2 (2 \| f_0^{\otimes k} \otimes f_1\|_q^2
        + 2 \langle f_0^{\otimes k} \otimes f_1, f_1 \otimes
        f_0^{\otimes k} \rangle_q)\\
&=
        2 \sum_{k=1}^\infty w_{k 1}^2 (1 + q^k) [k]_q!\\
&=
        2 + 2 \sum_{k=1}^\infty w_{k 1}^2 q^k [k]_q!.
\end{align*}
To prove the theorem it remains to show that 
\[
\sum_{k=1}^\infty w_{k 1}^2 q^k [k]_q! < q\qquad \text{for $q \in (0, 1)$.} 
\]
To this end note that $q^k < q$ for $k \geq 2$ and $q \in (0, 1)$, so:
\[
\sum_{k=1}^\infty w_{k 1}^2 (q^k - q) [k]_q! < 0
\]
from which it follows that:
\[
\sum_{k=1}^\infty w_{k 1}^2 q^k [k]_q! < q \sum_{k=1}^\infty w_{k 1}^2 [k]_q!
= q \| \widetilde{W} e_1\|_q^2 = q.
\]
We conclude that $\langle \Omega, (\gamma(X_0)+X_1)^4 \Omega\rangle_q
< \langle \Omega, (X_0+X_1)^4 \Omega\rangle_q$ for $q \in (0,1)$.
\end{proof}

The content of theorem~\ref{hoofdstelling} is shown graphically in
figure~\ref{mom_figuur} where the fourth moment of $X_0+X_1$ and a
numerical approximation of the fourth moment of $\gamma(X_0)+X_1$ are
plotted.

\end{document}